\documentclass[12pt]{article}
\usepackage[utf8]{inputenc}
\usepackage{color}
\usepackage{mathrsfs}
\usepackage{pifont}
\usepackage{amsmath}
\usepackage{amsthm}
\usepackage{txfonts}
 \usepackage{geometry}
\usepackage{latexsym}
\usepackage{amssymb}
\usepackage{graphicx}%
\usepackage{ifthen}
\usepackage{fancyhdr}
\usepackage{afterpage}
\usepackage{xcolor}
\usepackage{ifthen}

\newlength{\defbaselineskip}
\newcommand{\setlinespacing}[2]%
           {\setlength{\baselineskip}{#1 \defbaselineskip}}
\newcommand{\doublespacing}{\setlength{\baselineskip}%
                           {1.5 \defbaselineskip}}

\begin{document}
\doublespacing
\begin{center}
{\Large {\bf Sequential Detection of Transient Signals with Exponential Family Distribution }} \\
{\it Yanhong Wu{\footnote{Address correspondence to Yanhong Wu, Department of Mathematics, California State University Stanislaus, Turlock, CA 95382, USA; E-mail: ywu1@csustan.edu}}}\\
Department of Mathematics, California State University Stanislaus,\\
Turlock, California, USA \\
\end{center}

\begin{center}
ABSTRACT
\end{center}
\begin{quote}
We first consider the sequential detection of transient signals by generalizing the moving average chart to exponential family and study the false detection probability (FDP) and power of detection (POD) in the steady state. Then windowed adjusted signed (or modified directed) likelihood ratio chart is studied by treating it as normal random variable.
In the multi-parameter exponential family, the detection of the transient change of one of the  canonical parameters or a function of canonical parameters is considered by using  the generalized adjusted signed likelihood ratio chart. Comparisons with window restricted CUSUM and Shiryayev-Roberts (S-R) procedures show that the generalized signed likelihood ratio chart performs quite well.  Several important examples including the mean or variance change under normal model and a real example  are used for illustration.  \\
\end{quote}

\noindent{\bf Key Words:} Exponential family; Adjusted signed likelihood ratio; False detecting probability; Power of Detection; CUSUM and S-R procedures.

\bigskip

\section{Introduction}
Motivated from sequential detection of transient signals in a data stream, we consider the following transient change point model in a standard exponential family. Assume $\{X_1, .., X_t,...$ are independent random variables that follow the density function $f_{\theta} (x)=f_0(x) \exp(\theta x -c(\theta))$ for $i=\nu+1, ..., \nu+L$ and $f_0(x)$ for other $i$'s, where $f_0(x)$ is the baseline density function and $c(\theta)$ is differentiable for $|\theta|\leq K >0$ with $c(0) =c'(0)=0$. Our goal is to detect the transient signal sequentially. In quality control, the moving average (MA) or equivalently, the moving sums (MOSUM)  have been used for detecting a change (Bauer and Hackl (1978, 1980), Lai (1974),1995)) as an alternate to CUSUM and EWMA charts. Applications to test of parameter constancy in regression models are considered in Chu, et al. (1995). Recently, Noonan and Zhigljavsky (2020) studied the power of the MOSUM test under the normal model and Xie, et al. (2021) also studied a MA chart based on the largest eigenvalue of moving covariance matrix for detecting independent random mean change. 

In this communication, we first consider the performance of MA chart in the exponential family by studying its false detection probability (FDP) and power of detection (POD) defined under the following sense. 

Define $\bar{X}_{t;w} =\frac{1}{w}\sum_{t-w+1}^t X_i $ as the moving average of window width $w$. The MA chart makes an alarm at 
\[
\tau_{MA} = \inf\{t>0: \bar{X}_{t;w} >h\},
\]
where $h$ is typically the reference value for the mean of signal strength $c'(\theta)$ for some $\theta>0$.  Suppose when the change occurs, $\bar{X}_{t;w} $ is at the stationary state. We denote $P_{\theta}^*(.)$ as this corresponding probability measure when the parameter changes to $\theta>0$. Then we evaluate the performance of the MA chart by its FDP defined as $P_0^*(\tau_{MA} \leq L)$ and POD defined as $P_{\theta}^*(\tau_{MA} \leq L)$. 

In Section 2, we first give the approximation for FDP and POD for the MA chart in the exponential family. Our main contribution is to generalize the MA chart to the the situation when the signal length or signal strength are unknown. In Section 3, we first consider the generalized MA chart where a range for the signal length is assumed known. The window restricted adjusted signed likelihood ratio chart is used as a monitoring chart without assuming a reference value for the signal strength. The corresponding  FDP and POD are studied. By using the second order normal approximation for the adjusted signed likelihood ratio statistic, we can obtain accurate approximation for its FDP with the help of the localization technique after a likelihood ratio transformation (See Siegmund, et al. (2010, 2011)). The comparisons of the POD along with the window restricted CUSUM and S-R charts are conducted by using examples for detecting rate decrease in the exponential distribution and variance increase in the normal distribution. In Section 4, detection of a change in one or a function of the canonical parameters in a multi-parameter exponential family model is considered by using the general adjusted signed likelihood ratio statistic. Examples for mean increase or variance increase in the normal distribution are used to compare the PODs between the several charts. The results are extended further to multi-parameter case by using the chi-square approximation after a Bartlett adjustment. Detecting a change in both mean and variance in the normal distribution is used as an example. Real data from IBM stock price are used for demonstration. 

\section{FDP and POD of MA Chart}
By writing $Z_{t;w} =\sqrt{w}\bar{X}_{t;w}$ and $b=h\sqrt{w}$, we can rewrite
\[
\tau_{MA} = \inf\{t>0: Z_{t;w}>b\}
\]
The following the theorem gives the approximation for FDP. 

\noindent{\bf Theorem 1.} {\it  As $b, w \rightarrow \infty$ such that $b/\sqrt{w} \rightarrow h>0$ and $L=O(w)$, 
\begin{equation}
P_0^*(\tau_{MA} \leq L )
\approx \frac{Lb e^{-\theta b +wc(\theta/\sqrt{w})}}{\sqrt{2\pi} w (c''(\theta/\sqrt{w}))^{1/2}}  e^{-\frac{\theta}{\sqrt{w}} \rho_+}, 
\end{equation}
where $\rho_+ =E_0S_{\tau_+}^2/(2E_0S_{\tau_+})$ is the mean overshoot with
\[
\tau_+ =\inf\{ n>0: S_n = \sum_{i=1}^n (X_i-X_i')>0\}
\]
and $X_i'$ are i.i.d. copies of $X_i$ for $i=1,2,...$ under $P_0(.)$.
}

\noindent{\em Proof.} We use the following change of measure technique. For $1\leq t \leq L$,
define a changed measure $P_t(.)$ that under $P_t(.)$, $X_{t-w+1}, ...,X_t$ follow $f_{\theta/\sqrt{w}}(x)$ and other $X_i$'s follow $f_0(x)$. The log-likelihood ratio of $X_1, ..., X_L$ under $P_t(.)$ with respect to $P_0^*(.)$ is
\[
l_t = \theta Z_{t;w} -wc(\theta/\sqrt{w}).
\]
Note that $E_t(Z_{t;w}) = \sqrt{w}c'(\theta/\sqrt{w})$. So the value of $\theta$ will be selected such that the boundary $b$ equals to the mean of $E_t(Z_{t;w})$, i.e. $b= \sqrt{w} c'(\theta /\sqrt{w})$. Due to the monotone property of $c'(\theta)$, we can assume the solution exists and is unique. Denote $\tilde{l}_t= \theta ( Z_{t;w} -b)$. The change of measure technique in Siegmund and Yakir (2000) shows that
\begin{eqnarray*}
P_0^*(\tau_{MA} \leq L ) &=& P_0^*(\max_{1\leq t\leq L} Z_{t;w} \geq b ) \\
&=& \sum_{t=1}^L E_t [\frac{1}{\sum_1^L \exp(l_s^w)};\max_{1\leq t\leq L} Z_{t;w} \geq b  ] \\
&=& e^{\tilde{l}_t -l_t} \sum_{t=1}^L E_t [\frac{M_t}{S_t} e^{-(\tilde{l}_t +\ln(M_t)}; \tilde{l}_t +\ln(M_t) \geq 0 ],
\end{eqnarray*}
where $M_t =\max_{1\leq s\leq L} e^{l_s -l_t}$ and $S_t =\sum_{1\leq s\leq L} e^{l_s -l_t}$.

As $w \rightarrow \infty$, $\tilde{l}_t $ is asymptotically normal with mean 0 and variance $\theta^2 c''(\theta/\sqrt{w})$ under $P_t(.)$. By applying the localization technique in Siegmund, et al. (2010) along with the saddle-point approximation, we have 
\[
P_0^*(\tau_{MA} \leq L ) \approx L\frac{e^{-\theta b +wc(\theta/\sqrt{w})}}{\sqrt{2\pi} \theta (c''(\theta/\sqrt{w}))^{1/2}} \lim_{L\rightarrow \infty} E_t[\frac{M_t}{S_t}].
\]
In the following, we study in detail the approximation for $ E_t[\frac{M_t}{S_t}]$ in the exponential family case. For $s<t$ and $s \geq t-w$,
\[
l_s-l_t =\frac{\theta}{\sqrt{w}} (\sum_{s-w+1}^{t-w} X_i -\sum_{s+1}^t X_i) =\frac{\theta}{\sqrt{w}}
\sum_{s+1}^t (X_{i-w} -X_i) ,
\]
which is a two-sided random walk with drift
\[
E_t[l_s-l_t]= - \frac{\theta}{\sqrt{w}} (t-s) c'(\frac{\theta}{\sqrt{w}}) =-\frac{\theta b}{w} |t-s|,
\]
and variance
\[
Var_t(l_s-l_t)= |t-s| \frac{\theta^2}{w} (c''(0)+c''(\frac{\theta}{\sqrt{w}} )).
\]
Same result holds for $s>t, s\leq t+w$. 
Thus we have the following approximation
\[
E_t[\frac{M_t}{S_t}] \approx \frac{b \theta}{w} \upsilon(\frac{\theta}{\sqrt{w}} (c''(0)+c''(\theta/\sqrt{w}))^{1/2} ),
\]
where the last factor equals to $\lim_{d \rightarrow \infty} E\exp(-(\theta/\sqrt{w})(S_{\tau_d}-d))$ with 
\[
\tau_d =\inf\{ t>0: S_t =\sum_{k=1}^t (X_i-X_i') >d\} 
\]
where $X_i$ follows $f_{\theta} (x)$ and $X_i'$ follows $f_0(x)$ and can be approximated as
\[
\upsilon(\frac{\theta}{\sqrt{w}} (c''(0)+c''(\theta/\sqrt{w}))^{1/2} ) 
\approx e^{-\frac{\theta}{\sqrt{w}} \rho_+}.
\]
Theorem 1 is proved by combining the above approximations,

Note that typically $h$ is quite small, thus $\theta/\sqrt{w} $ is also small. The above approximation can further be simplified as 
\begin{equation}
P_0^*(\tau_{MA} \leq L ) \approx \frac{Lb}{w\sqrt{c''(0)}}\phi(\frac{b}{\sqrt{c''(0)}} ) e^{-b\rho_+/(\sqrt{w} c''(0))} = \frac{Lh}{\sqrt{wc''(0)}}\phi(\frac{h\sqrt{w}}{\sqrt{c''(0)}} ) e^{-h\rho_+/c''(0)},
\end{equation}
where $\phi(x)$ is the standard normal density function.  \qed

\noindent{\bf Remark.} In the normal case, to detect a change of mean from 0 to a positive value, $\rho_+ \approx 0.5826\sqrt{2} \approx 0.824$ and 
\begin{equation}
P_0^*(\tau_{MA} \leq L ) \approx \frac{Lb}{w} \phi (b ) e^{ - 0.5826\sqrt{2}b/\sqrt{w}} \approx \frac{Lh}{\sqrt{w}} \phi (h \sqrt{w} ) e^{ - 0.824 h },
\end{equation}
as shown in Wu and Siegmund (2022).

\noindent{\bf Example 1.} (Rate change) Suppose the original observations $Y_i$ follow $Exp(\lambda_0)$ for $i\leq \nu$ and $Exp(\lambda )$ for $i > \nu$ and $\lambda <\lambda_0$. By making a transformation $X_i = \lambda_0 Y_i -1$ and defining $f_0(x) =\exp(-(x+1))$ for $x \geq -1$, we see that for $\theta = 1-\lambda /\lambda_0$,
\[
c(\theta)= -\theta -\ln (1-\theta) =-(1-\lambda /\lambda_0) -\ln (\lambda/\lambda_0).
\]
This gives $c'(\theta )= \theta/(1-\theta)$  and $c''(\theta)=1/(1-\theta)^2$. From $c'(\theta/\sqrt{w})=b/\sqrt{w}$, we have
\[
\theta =b/(1+b\sqrt{w}).
\]
Equation (2.1) gives 
\begin{equation}
P_0^*(\tau_{MA} \leq L )\approx \frac{Lb}{\sqrt{2\pi} w} (1+ \frac{b}{\sqrt{w}})^{w-1} e^{-\sqrt{w} (b+\rho_+/w)}.
\end{equation}
In the exponential case, because of the memoryless property, $\rho_+ =1$. 

\noindent{\bf Example 2.} (Variance change) Suppose the original data $Y_i$ follow $N(0,\sigma_0^2)$ for $i\leq \nu$ and  $N(0, \sigma^2) $ for $i > \nu$ and $\sigma^2 >\sigma_0^2$. By defining $X_i = (Y_i^2/\sigma_0^2 -1)/\sqrt{2}$ and $f_0(x)$ as the density function under $N(0, \sigma_0^2)$, then for $\theta = (1-\sigma_0^2/\sigma^2)/\sqrt{2}$, we have
\[
c(\theta)= -\frac{\theta}{\sqrt{2}} -\frac{1}{2} \ln (1-\sqrt{2}\theta) =-\frac{1}{2} (1-\frac{\sigma_0^2}{\sigma^2}) -\frac{1}{2} \ln (\sigma_0^2/\sigma^2).
\]
Note $c'(\theta)= \theta /(1-\sqrt{2}\theta)$ and $c''(\theta) =1/(1-\sqrt{2}\theta)^2$. From
$c'(\theta/\sqrt{w})=b/\sqrt{w}$, we get
\[
\theta = b/(1+\sqrt{2}b/\sqrt{w}).
\]
Equation (2.1) gives
\begin{equation}
P_0^*(\tau_{MA} \leq L )\approx \frac{Lb}{\sqrt{2\pi} w} (1+\sqrt{2} \frac{b}{\sqrt{w}})^{w/2-1} e^{-\sqrt{w} (b+\sqrt{2}\rho_+/w)/\sqrt{2}}.
\end{equation}
For the variance change, $\rho_+ \approx 1.167$ by simulation.

Next, we study the POD for the MA chart by separating two cases. When $c'(\delta) <h$, let $w \rightarrow \infty$ and $L\rightarrow \infty$ at least at the same order as $w$. By treating $1/w$ as the time scale, the local POD can be approximated as
\begin{eqnarray*}
P_{\delta}^*(\tau_{MA} \leq L ) = && P_0^*( \max_{1\leq t \leq L} (\bar{X}_{t;w} + \min(t/w,1) c'(\delta )) >h ) \nonumber \\
= &&    \int_0^{L/w} \frac{\sqrt{w} (h - \min(u, 1) c'(\delta))}{\sqrt{c''(0)}} \phi( \frac{\sqrt{w} (h - \min(u, 1) c'(\delta))}{\sqrt{c''(0)}} ) e^{-\rho_+  (h - \min(u, 1) c'(\delta))/c''(0)} du  
\end{eqnarray*}

When $c'(\delta) \geq h$, by ignoring the overshoot at the crossing time $t$, we have
\[
\bar{X}_{t;w} +\min(\frac{t}{w}, 1) c'(\delta) =h .
\]
This gives
\[
\frac{t}{w} =\frac{h}{c'(\delta)} -\frac{\bar{X}_{t;w}}{c'(\delta)}.
\]
As $w\rightarrow \infty$,
\[
E_{\delta}^*(\tau_{MA}) \approx \frac{h}{c'(\delta)} w, ~~~
Var_{\delta}^*(\tau_{MA}) = \frac{wc''(\delta)}{(c'(\delta))^2} .
\]
The normal approximation gives the following approximation:
\[
P_{\delta}^*(\tau_{MA} \leq L) \approx \Phi \left( \frac{c'(\delta)L-hw}{\sqrt{wc''(\delta)}} \right).
\]
A more accurate approximation can be obtained by using a continuous correction of $h$ from the approximation (2.2). 

\section{Generalizations}
There are several shortcomings of the MA chart. One is that the signal length is typically unknown. The other is that the design may be complicated for different models under the exponential family. 
In this section, we consider two possible generalizations of the simple MA chart. One is to generalized MA chart by covering an interval $[w_0, w_1]$ for possible signal length. The second is to use a second order corrected signed likelihood ratio statistic to transform the monitoring statistic into normal. 

\subsection{Generalized MA chart}

When the signal length $L$ can be estimated within an interval $[w_0,w_1]$, we can use the following GMA chart by making an alarm at
\[
\tau_{GMA} =\inf\{t>0: \max_{w_0\leq w\leq w_1} Z_{t;w} \geq b\}
\]

By using the same technique as in Section 2, we have the following approximation for the FDP.
\[
P_0^*(\tau_{GMA}\leq L ) \approx L \int_{w_0}^{w_1} \frac{\theta b^2}{4w^2} \frac{e^{-\theta b +wc(\theta/\sqrt{w})}}{ (c''(\theta/\sqrt{w}))^{1/2}}  e^{-\frac{\theta}{\sqrt{2 w}} \rho_+} dw ,
\]
where $\theta $ inside the integration satisfying $b/\sqrt{w} =c'(\theta/\sqrt{w})$. For small $\theta/\sqrt{w_0}$, the above approximation can be simplified as 
\[
P_0^*(\tau_{GMA}\leq L ) \approx \frac{Lb}{\sqrt{c''(0)}} \phi(\frac{b}{\sqrt{c''(0)}}) \int_{b/\sqrt{w_1c''(0)}}^{b/\sqrt{w_0c''(0)}} \frac{u}{2} e^{-2\rho_+ u} du .
\]
A proof can be obtained as in Siegmund, et al. (2011) or Wu and Siegmund (2022).  

\subsection{Adjusted Signed Likelihood Ratio (LR) Chart}

Under the exponential family, we can match the adjusted signed likelihood ratio (LR) statistic to a standard normal random variable up to the second order. 

For observations $X_{t-w+1}, ..., X_t$ that follow $f_0(x)$, the signed (or directed) likelihood ratio is defined as
\[
R_t = sgn (\hat{\theta}_t) (2 l(\hat{\theta}_t))^{1/2} =
sgn (\hat{\theta}_t) (2 w (\hat{\theta}_t \bar{X}_{t;w} -c(\hat{\theta}_t)))^{1/2},
\]
where $\hat{\theta}_t$ satisfies $c'(\hat{\theta}_t) =\bar{X}_{t;w}$. 
By defining $U_t= \hat{\theta}_t\sqrt{w c''(\hat{\theta})}$ , the adjusted signed (or modified directed) LR statistic is defined by
\[
R_t^* =R_t +\frac{1}{R_t} \ln (U_t/R_t). 
\]
Note that $R_t^*$ is a function of $\bar{X}_{t;w}$ and $P_0(R_t^* \leq t) =\Phi(t)$ with order of $O(w^{-3/2})$ for $t=O(1)$ and $O(1/w)$ for $t=O(\sqrt{w})$ (e.g. Jensen (1995, pg 118, (5.2.2)). We shall treat $R_t^*$ as a standard normal variable and define the LR chart by making an alarm at
\[
\tau_{LR} =\inf\{ t>0: R_t^* > b\}.
\]
An approximation for FDP is given in the following theorem.

\noindent{\bf Theorem 2.} {\it As $w, b \rightarrow \infty $ such that $b/\sqrt{w} \rightarrow h>0 $ and $L=O(\sqrt{w})$,
\begin{equation}
P_0^*(\tau_{LR} \leq L ) \approx \frac{L\theta c'(\theta)}{b} \phi(b) \upsilon (\theta  \sqrt{c''(0)+c''(\theta)} )  \approx  \frac{Lb}{w}\phi(b ) e^{-b\rho_+/\sqrt{w c''(0)}},
\end{equation}
where $\rho_+ $ is defined as in Theorem 1. 
}

\noindent{\em Proof.} 
We extend the techniques as in Section 2. Define a changed measure $P_t(.)$ such that $X_{t-w+1},..., X_t$ follow $f_{\theta}(x)$ and others follow $f_0(x)$ such that under $P_t(.)$, the likelihood ratio can be approximated by
\[
\frac{dP_t}{dP_0}(R_t^*) \approx \exp(b R_t^* -b^2/2 ) (1+O(w^{-3/2})).
\]
The value of $\theta$ will be decided by equating the mean $E_t^* R_t^* =b$, the boundary. Thus, the log-likelihood ratio is $l_t =b R_t^* -b^2/2$. 
Denote $\tilde{l}_t = b(R_t^* -b)$ with mean 0 and variance $b^2$. 
By applying the localization technique, for $L =O(w)$, we can immediately get the following approximation
\begin{eqnarray*}
P_0^*(\tau_{LR} \leq L) &=& P_0^*(\max_{1\leq t\leq L} R_t^* \geq b ) \\
&=& \sum_{t=1}^L E_t [\frac{1}{\sum_1^L \exp(l_s^w)};\max_{1\leq t\leq L} R_t^* \geq b  ] \\
&=& e^{\tilde{l}_t -l_t} \sum_{t=1}^L E_t [\frac{M_t}{S_t} e^{-(\tilde{l}_t +\ln(M_t)}; \tilde{l}_t +\ln(M_t) \geq 0 ] \\
&\approx &  \frac{L}{b} \phi (b) \lim_{L\rightarrow \infty} E_t [\frac{M_t}{S_t}].
\end{eqnarray*}
where $M_t=\max_{0\leq s \leq L} \exp((b(R_s^* -R_t^*))$ and $S_t =\sum_{0\leq s\leq t}\exp((b(R_s^* -R_t^*)) $. 

As $w\rightarrow \infty$, since $E_t R_t^* \approx b$, by assuming $b/\sqrt{w} \rightarrow h $, by using the stochastic expansion under $P_t(.)$ below, $\theta$ satisfies
\begin{equation}
\psi (\theta) =\theta c'(\theta ) -c(\theta)  \approx b^2/2w .
\end{equation}
When $b^2/w$ is small, $\theta \approx b/\sqrt{c''(0) w}$. To study  the limit of $E_t [\frac{M_t}{S_t}]$, we use the following stochastic expansion under the changed measure
$P_t(.)$ with parameter $\theta$
\begin{eqnarray*}
R_t &\approx & \approx (2w (\hat{\theta}_t \bar{X}_{t;w} -c( \hat{\theta}_t ))^{1/2} \\
& =& \sqrt{2w} (\psi (\theta) +(\hat{\theta}_t -\theta) c'(\theta) +\hat{\theta}_t (\bar{X}_{t;w} -c'(\theta)) -(c(\hat{\theta}_t) -c(\theta)))^{1/2} \\
& =& \sqrt{2w} (\psi (\theta) +\theta(\bar{X}_{t;w} -c'(\theta))+ (\hat{\theta}_t -\theta)(c'(\hat{\theta}_t ) -c'(\theta)) -\frac{(\hat{\theta}_t -\theta))^2}{2} c''(\theta) +o_p(1/w ))^{1/2} \\
&=& \sqrt{2w} (\psi (\theta) +\theta(\bar{X}_{t;w} -c'(\theta)) +\frac{(\hat{\theta}_t -\theta))^2}{2} c''(\theta) +o_p(1/w ))^{1/2} \\
&=& \sqrt{2w \psi(\theta)} \left( 1+ \frac{\theta}{\psi (\theta)} (\bar{X}_{t;w} -c'(\theta)) +\frac{c''(\theta)}{2\psi (\theta)} (\hat{\theta}_t -\theta))^2 +o_p(1/w )\right)^{1/2} \\
&=& \sqrt{2w \psi(\theta)} \left( 1+ \frac{\theta}{2\psi (\theta)} (\bar{X}_{t;w} -c'(\theta))
+\frac{c''(\theta)}{4\psi (\theta)} (\hat{\theta}_t -\theta))^2 -\frac{\theta^2}{8\psi^2(\theta)} (\bar{X}_{t;w} -c'(\theta))^2 +o_p(1/w) \right ) \\
&=& \sqrt{2w \psi(\theta)} +\sqrt{\frac{w}{2\psi (\theta)}} \theta  (\bar{X}_{t;w} -c'(\theta)) +o_p(1) \\
& =& b+ \frac{w}{b} \theta  (\bar{X}_{t;w} -c'(\theta)) +o_p(1) . 
\end{eqnarray*}
When $\theta $ is small,
\[
R_t \approx b+ \frac{\sqrt{w}}{\sqrt{c''(0))}} \bar{X}_{t;w} (1+o_p(1)). 
\]
For $s= t \pm O(\sqrt{w})$, since $E_t \bar{X}_{s;w} = (1-\frac{|t-s|}{w})c'(\theta) $, we see the same expansion still holds for $R_s$.
On the other hand, at the first order, $E_tR_t \approx b =O(\sqrt{w})$ and 
$E_t U_t = \theta \sqrt{w c''(\theta)} (1+o(1)) $. Thus
\[
E_t [\frac{1}{R_t} \ln \frac{U_t}{R_t} ] = \frac{1}{b} \ln \frac{ \theta \sqrt{w c''(\theta)}}{b} (1+o(1)).
\]

Thus, just like the MA chart, as $w\rightarrow \infty $, 
\[
R_s^* -R_t^* \approx R_s -R_t  \approx \frac{w}{b} \theta (\bar{X}_{s;w} -\bar{X}_{t;w}).
\]
It behaves like a two-sided random walk with drift 
\[
E_t[l_s-l_t] \approx - \theta |t-s)| c'(\theta) \approx -\frac{b^2}{wc''(0)} |t-s|,
\]
and variance
\[
Var_t(l_s-l_t)\approx  \theta^2 |t-s| (c''(0)+c''(\theta)) \approx 2|t-s| \frac{b^2}{wc''(0)}.
\]
The theorem is proved by using the above approximations. \qed 

To study the POD, for any signal $\delta>0$, we derive the mean and variance of $\tau_{LR}$ under the assumption that $b^2/(2w) \rightarrow b^{*2}>0 $ and
\[
\delta c'(\delta) -c(\delta ) \geq  b^{*2}.
\]
At the boundary crossing time  $t$, $R_t^* \approx R_t \approx b$. This gives
\[
\hat{\theta}_t (\min(t/w, 1) c'(\delta) +\bar{X}_{t;w}) c(\hat{\theta}_t) = b^{*2}.
\]
Denote $\delta^* \leq \delta $ such that $\delta^* c'(\delta^*) -c(\delta^*) =b^{*2}$ and
$t^* =c'(\delta^*)/c'(\delta)$. Then as $\bar{X}_{t;w} \rightarrow 0$, we can see that $ t/w \rightarrow t^*$ and $\hat{\delta}_t \rightarrow \delta^*$.

By writing 
\[
\frac{t}{w} c'(\delta) = t^* c'(\delta) +(\frac{t}{w} -t^*) c'(\delta),
\]
we have
\begin{eqnarray*}
\hat{\theta}_t &= &c'^{-1}(t^* c'(\delta)) +(\frac{t}{w} -t^*) c'(\delta)+\bar{X}_{t;w})  \\
&=&c'^{-1}(t^* c'(\delta)) +\frac{d}{dx}c'^{-1}(t^* c'(\delta))[(\frac{t}{w} -t^*) c'(\delta)+\bar{X}_{t;w})] +O_p(1/w),
\end{eqnarray*}
and
\begin{eqnarray*}
c(\hat{\theta}_t) &=&c(c'^{-1}(t^* c'(\delta)) +(\frac{t}{w} -t^*) c'(\delta)+\bar{X}_{t;w}))  \\
&=&c(c'^{-1}(t^* c'(\delta))) + t^* c'(\delta) \frac{d}{dx}c'^{-1}(t^* c'(\delta)) [(\frac{t}{w} -t^*) c'(\delta)+\bar{X}_{t;w})]+O_p(1/w).
\end{eqnarray*}
Combining the above expansions, we get
\[
\hat{\theta}_t (\frac{t}{w} c'(\delta) +\bar{X}_{t;w}) c(\hat{\theta}_t)
=\delta^*c'(\delta^*)-c(\delta^*) +c'^{-1}(t^* c'(\delta))[(\frac{t}{w} -t^*) c'(\delta)+\bar{X}_{t;w})]+O_p(1/w) .
\]
This implies
\[
(\frac{t}{w} -t^*) c'(\delta)+\bar{X}_{t;w}= O_p(1/w).
\]
From this, we obtain
\[
E_{\delta}^*[\tau_{LR}/w] \approx t^*; ~~~~ Var_{\delta}^* [\tau_{LR}/w] \approx \frac{c''(\delta)}{c'^2(\delta)}.
\]
The normal approximation can be used to get the first order results for POD. 
A more careful second order stochastic expansion can help to derive more accurate approximations. 

{\bf Remark.} The above adjusted signed likelihood ratio chart, in principle,  can be applied to any parametric model (e.g. Brazzale, et al. (2006, pg 17) or Severini (2000, pg 242)). Denote $l_{t;w} (\theta; \hat{\theta}_t, a) =\sum_{k=t-w+1}^t \ln f_{\theta}(X_k) $ as the log-likelihood where $(\hat{\theta}_t, a)$ is a transform of the data and $a_t$ is an asymptotic ancillary with $\hat{\theta}_t$ the MLE from $W_{t-w+1}, ..., X_t$. The signed log-likelihood ratio can be written as 
\[
R_t =sgn(\hat{\theta}_t -\theta_0) [2(l_{t;w}(\hat{\theta}_t) - l_{t;w}(\theta_0)]^{1/2}. 
\]
Then the adjusted signed log-likelihood ratio $R_{t;w}^*$ can be calculated with 
\[
U_t = \hat{j}^{-1/2} (l_{;\hat{\theta}_t} (\hat{\theta}_t) -l_{;\hat{\theta}_t} (\theta_0)),
\]
where $\hat{j} = -\frac{\partial^2}{\partial \theta^2} l_{t;w} (\theta; \hat{\theta}_t, a_t) |_{\theta=\hat{\theta}_t}$ and $l_{;\hat{\theta}_t} (\theta)  =\frac{\partial} {\partial \hat{\theta}_t} l_{t;w} (\theta; \hat{\theta}_t, a_t)$.

\subsection{Windowed CUSUM and S-R Procedures}

For detecting a change, the CUSUM and Shiryayev-Roberts (S-R) procedures have certain optimalities. 
In Yakir (1995) and Siegmund and Yakir (2008), a new approach to calculate the average run length of S-R procedure for detecting a change in an exponential family is used, where the FDP is obtained as a by-product. Define
\[
l_k^{\theta}(n)= \sum_{n-k+1}^n (\theta X_i -c(\theta)),
\]
and $l_k(n) =\ln (\int \exp(l_j^{\theta} (n)) \rho (\theta) d\theta ) $ for a prior distribution $\rho(\theta)$ over $ |\theta|\leq K$. The window restricted S-R procedure makes an alarm at
\[
\tau_{SR}^w =\inf\{ n>0: R_{n;w} = \ln \sum_{k=1}^w e^{l_k(n)} >c \}.
\]
The Equation (8) in Siegmund and Yakir (2008) shows that for $c < < w$ and $c^{1+\epsilon} < L < \exp (c^{\epsilon} )$
\[
P_0^*(\tau_{SR}^w \leq L ) \approx Le^{-c} \int \lambda(\theta) \rho(\theta ) d\theta,
\]
where $\lambda (\theta )= \lim_{a\rightarrow \infty} E_1^{\theta} (\exp(-(l_1^{\theta}(\tau^{\theta}) -a))$ with $\tau^{\theta} =\inf\{n: l_1^{\theta}(n)\geq a\}$.

The same method can  be used to obtain the similar results based on CUSUM procedure, see Siegmund (1988) in the normal case. In practice, we assign a reference value for the changed parameter $\theta$, say $\delta>0$. In this situation, the windowed S-R and CUSUM procedures are defined as
\[
\tau_{SR}^w =\inf\{t>w: R_{t;w}(\delta) = \ln \sum_{k=1}^w \exp(l_k^{\delta} (t)) >c\}
\]
\[
\tau_{CS}^w =\inf\{t>w: \max_{1\leq k \leq w} l_k^{\delta} (t) >d\}
\]
Approximations for FDP can be obtained as in Yakir and Siegmund (2008). For example, for the S-R process,
\[
P_0^*(\tau_{SR}^w \leq L) \approx L e^{-c} \lambda(\delta ),
\]
and for the CUSUM procedure, it is given by
\[
P_0^*(\tau_{CS}^w \leq L) \approx L(\delta c'(\delta) -c(\delta)) e^{-d} \lambda^2(\delta).
\]

\noindent{\bf Example 1.} (Cont'd)  For the exponential case, 
\[
R_t =\sqrt{2w} sgn(\bar{X}_{t;w}) (\bar{X}_{t;w} -\ln (1+\bar{X}_{t;w}))^{1/2},
\]
and 
\[
U_t= \sqrt{w} \bar{X}_{t;w}.
\]
Therefore,
\[
R_t^*= R_t+\frac{1}{2R_t} \ln \frac{(\bar{X}_{t;w})^2}{2 (\bar{X}_{t;w} -\ln (1+\bar{X}_{t;w}))}.
\]
On the other hand, for a reference rate $\lambda_1 <\lambda_0$, the CUSUM process becomes
\[
\max_{1\leq k \leq w} l_k^{\delta} (t) = \max_{1\leq k \leq w} [(\lambda_0 -\lambda_1) \sum_{t-k+1}^t Y_i +k\ln (\lambda_1 \lambda_0 )],
\]
and the S-R process  becomes
\[
R_{t;w} =\ln \sum_{1\leq k \leq w} \exp[(\lambda_0 -\lambda_1) \sum_{t-k+1}^t Y_i +k\ln (\lambda_1 \lambda_0 )].
\]

\noindent{\bf Example 2. } (Cont'd) For the variance change case,
\[
R_t = sgn(\bar{X}_{t;w}) \sqrt{w} ( \sqrt{2} \bar{X}_{t;w} -\ln (1+\sqrt{2} \bar{X}_{t;w}))^{1/2},
\]
and $U_t =\sqrt{w} \bar{X}_{t;w}$. Thus
\[
R_t^* = R_t+\frac{1}{2R_t} \ln \frac{(\bar{X}_{t;w})^2}{\sqrt{2} \bar{X}_{t;w} -\ln (1+\sqrt{2} \bar{X}_{t;w})}.
\]
For a reference value $\sigma_1^2 >\sigma_0^2$, the CUSUM process with window width $w$ is given by
\[
\max_{1\leq k \leq w} l_k^{\delta} (t) = \max_{1\leq k \leq w} [(1/\sigma_0^2 -1/\sigma_1^2) \sum_{t-k+1}^t Y_i^2 +k\ln (\sigma_0^2/\sigma_1^2 )]/2,
\]
and the S-R process is 
\[
R_{t;w} =\ln \sum_{1\leq k \leq w} exp([(1/\sigma_0^2 -1/\sigma_1^2) \sum_{t-k+1}^t Y_i^2 +k \ln (\sigma_0^2/\sigma_1^2 )]/2).
\]

Note that the MA chart is based on Rao's score statistic, $U_t$ is Wald's statistic, and $R_t$ is the signed likelihood ratio statistic. At the first order, all three are asymptotically normal. 

In Examples 1 and 2, the Rao's statistic is the same as $U_t$. To compare the PODs between the several charts, we conduct a simulation comparison of POD for same FDP =0.02 between the MA chart, the signed likelihood ratio statistic, and the adjusted signed likelihood ratio statistic for $w=20$ and $L=10, 20, 30 $. All simulations are replicated 10,000 times. We can see that all charts give roughly the same POD.  Table 1 lists the POD's in the exponential case for $\lambda =1$ as initial rate with data transformation $Y_i -1$ and $1/\lambda =1.25, 1.50, 1.75, 2.0, 2.25, 2.5, 2.75, 3$ for changed rates. Included are also the CUSUM and S-R procedures with reference value $\delta =0.5$ ($\lambda =0.5$). The numbers in the brackets are the approximated values of FDP based on approximations (2.2) and (3.7) with normal approximations. We see that MA chart has much heavier tail; while $R_t$ has light tail.

\begin{table}[ht]
\begin{center}
Table 1.  Comparison of POD for $w=20$ for rate decrease in exponential R.V. \\
\begin{tabular}{|c| c |c |c |c |c|c|}   \hline
L &$1/\lambda$ & MA & $R_t$ & $R_t^*$ & CUSUM & S-R\\
 & boundary & b=3.10 & b=2.55 & b=2.67 &d=4.48 & c=6.14\\ \hline
20 & 1.00 & 0.0199(0.006) & 0.0191(0.0246) & 0.0200(0.0185)& 0.0196 &0.0200\\
 &1.25 & 0.1127 & 0.1224 & 0.1112&0.1174&0.1296 \\
& 1.50 & 0.3359 & 0.3465 & 0.3274& 0.3380 &0.3481\\
& 1.75 & 0.5864 & 0.6038 & 0.5849& 0.5899 &0.5962\\
& 2.00 & 0.7780 & 0.7949 & 0.7826&0.7786&0.7720 \\
& 2.25 & 0.8993 & 0.8966 & 0.8994 &0.8888  &0.8931 \\
& 2.50 & 0.9539 & 0.9580 & 0.9516 &0.9496 & 0.9481\\
& 2.75 & 0.9774 & 0.9811 & 0.9788 &0.9768 & 0.9784 \\
& 3.0 & 0.9895 & 0.9916 & 0.9904 &0.9886 &0.9905 \\ \hline
30 &1.0 & 0.0302 & 0.0298&0.0282 &0.0282 &0.0289 \\
  & 1.25 & 0.1952 & 0.2111 & 0.1923 &0.1825 &0.1926 \\
  & 1.50 & 0.5136& 0.5338 & 0.5161 &0.4901 &0.5015 \\
  & 1.75 & 0.7929 & 0.8000 &0.7843 &0.7560 & 0.7659 \\
  &2.00 & 0.9249  & 0.9324 &0. 9227 &0.9070 & 0.9122\\
  &2.25 & 0.9780 & 0.9752  & 0.9759 & 0.9703 & 0.9691 \\
  & 2.50 & 0.9937&0.9954 & 0.9921 &0.9903  & 0.9926 \\
  &2.75 & 0.9991 & 0.9979 & 0.9987 &0.9965 &0.9976 \\
  & 3.00 & 0.9995 &0.9992 & 0.9996 &0.9989 &0.9994 \\  \hline 
10 &1.0 & 0.0124 &0.0134 &0.0126 & 0.0140 &0.0126 \\
   & 1.25 & 0.0403 &0.0412 &0.0341 &0.0509 &0.0540 \\
   & 1.50 & 0.0980 &0.0973 &0.0969 &0.1406 &0.1525 \\
   &1.75 & 0.1895 &0.1870  &0.1816 &0.2783 & 0.2792 \\
   &2.00 &0.2983 &0.3188  &0.3039 &0.4283  & 0.4398 \\
   &2.25 & 0.4370 &0.4469 &0.4245 &0.5632  &0.5914\\
   &2.50 & 0.5471 &0.5613 &0.5463 &0.6865  &0.6898 \\
   & 2.75 & 0.6517 &0.6592 &0.6498 &0.7818 & 0.7821 \\
   &3.00 & 0.7309 &0.7353 &0.7345 &0.8386  & 0.8426 \\ \hline 
\end{tabular}
\end{center}
\end{table}

In Table 2, we also conduct a simulation comparison of POD for detecting variance increase from 1 as in Example 2. The transformation is $(Y^2-1)/\sqrt{2}$. For the CUSUM and S-R procedures, we take the reference value $\delta= 0.5/\sqrt{2}$ that corresponds to $\sigma^2 =2$. We see that the chart based on $R_t$ underestimate the alarm limit  and the MA chart has heavier tail as the alarm limit is much larger.

\begin{table}[ht]
\begin{center}
Table 2.  Comparison of POD for $w=20$ for variance increase in normal R.V. \\
\begin{tabular}{|c |c |c |c |c |c |c|}   \hline
L & $\sigma^2$ & MA & $R_t$ & $R_t^*$ & CUSUM &S-R\\
  & boundary & b= 3.30 & b=2.55 & b=2.65& d=3.96 & c=5.84 \\ \hline
20 &1.00 & 0.0191(0.003) & 0.0189(0.0246) & 0.0188(0.0194) & 0.0192 &0.0190 \\
 & 1.25 & 0.0713 & 0.0737 & 0.0751&0.0797 & 0.0826 \\
 & 1.50 & 0.1812 & 0.1875 & 0.1846 &0.1946 &0.2070  \\
 &1.75 & 0.3362 & 0.3436 & 0.3403 &0.3617 &  0.3712\\
 & 2.00 & 0.4932 & 0.4995 & 0.4983 & 0.5120 & 0.5165 \\
 &2.25 & 0.6306 & 0.6355 & 0.6359 & 0.6517& 0.6545 \\
 &2.50 & 0.7395 & 0.7445 & 0.7436 &0.7449 &  0.7629  \\
 &2.75 & 0.8178 & 0.8239 & 0.8205 & 0.8251 & 0.8408 \\
 & 3.0 & 0.8739 & 0.8761 & 0.8767 &0.8758& 0.8802 \\ \hline
30 & 1.0 & 0.0284 & 0.0257 &0.0258 & 0.0296 &0.0264 \\
   & 1.25 & 0.1205 & 0.1290 &0.1176 & 0.1253& 0.1270 \\
   & 1.50 & 0.3053 & 0.3137 &0.3109 & 0.3108 & 0.3078 \\
   & 1.75 & 0.5171 & 0.5190 & 0.5162 & 0.5134 & 0.5205 \\
   & 2.00 & 0.6943 &0.6953  &0.6934 & 0.6892 & 0.6953 \\
   &2.25 & 0.8168  & 0.8199 & 0.8195 & 0.8154 & 0.8202 \\
   &2.50 & 0.8927 & 0.9039  &0.8955 &0.8927& 0.8915 \\
   & 2.75 & 0.9402 &0.9418 &0.9423 & 0.9370 & 0.9398 \\
   & 3.00 & 0.9712 & 0.9633 &0.9650 & 0.9685 & 0.9643 \\ \hline
 10 & 1.00 & 0.0109 &0.0126 &0.0120 &0.0126 & 0.0128 \\
    & 1.25 & 0.025 & 0.0253 &0.0270 &0.0375 & 0.0364 \\
    &1.50 & 0.0595 & 0.0560  &0.0593 & 0.0824 & 0.0896 \\
    &1.75 & 0.0976 & 0.1015 & 0.1075 & 0.1522 & 0.1571 \\
    &2.00 & 0.1676 & 0.1680 & 0.1707 & 0.2343 & 0.2440 \\
    &2.25 & 0.2274 & 0.2451 &0.2415  &0.3249  & 0.3336 \\
    &2.50 & 0.3093 & 0.3145 & 0.3137 &0.4055 & 0.4210 \\
    &2.75 & 0.3846 & 0.3920 & 0.3997  & 0.4899 & 0.5033 \\
    &3.00 & 0.4501 &0.4763  & 0.4591  & 0.5566 & 0.5665 \\ \hline 
\end{tabular}
\end{center}
\end{table}

\section{Detecting a Change of Scale Parameter in Multi-parameter Exponential Family Model}

For the multi-parameter exponential family, Siegmund and Yakir (2008) considered the minimax optimality of the mixture Shiryayev-Roberts procedure in the exponential family case in terms of average delay detection time. An approximation for FDP is given in Equation (8).  A comparison of variety of test charts including EWMA, MA, CUSUM, and GLRT in the normal case is conducted in Wu and Siegmund (2022). In the exponential family case, the detecting process will be formed by the windowed Rao-score statistic, Wald's statistic, and the likelihood ratio test statistic. However, a high order correction looks complicated. In the following, we shall focus on the detection of a scale parameter change by treating other parameters as nuisance. 

\subsection{Detection a scale canonical parameter change}
In this subsection, we first consider the detection of a scale canonical parameter change with other parameters being treated as nuisance and staying the same. 
We shall again use the adjusted signed likelihood ratio statistic by treating it as normal. 

Assume $\{X_1, .., X_t,...\}$ are independent random variables that follow the density function $f_{\theta} (x)=f_0(x) \exp(\theta t(x) -c(\theta))$ for $i=\nu+1, ..., \nu+L$ and $f_0(x)$ for other $i$'s, where $f_0(x)$ is the baseline density function, $\theta =(\theta_1, ..., \theta_p)=(\theta_1, \theta^{(2)})$, $t(x)=(t_1(x), ..., t_p(x))^T =(t_1(x), t^{(2)}(x))^T$ and $c(\theta)$ is differentiable for $||\theta||\leq K >0$ with $c(0)=0$ and $\partial c(0)/\partial \theta =(0, ....,0)$. 

We first consider the detection of change of $\theta_1$ from $\theta_{01} $ and $\theta^{(2)}$ is unknown and may  be randomly varying. For the observations  $X_{t-w+1}, ..., X_t$ with window width $w$, denote by
$\bar{t}_1 = (t_1(X_{t-w+1}) +...+t_1(X_t))/w$, $\bar{t}^{(2)} = (t^{(2)}(X_{t-w+1})+...+t^{(2)}(X_t))/w$. For convenience of notation, we denote 
$\frac{d}{d \theta^T} c(\hat{\theta}) =c'(\theta) $ and $\frac{\partial^2}{\partial\theta \partial \theta^T} c(\theta) =c''(\theta)$. Let
\[
c'(\hat{\theta}) =(\bar{t}_1, \bar{t}^{(2)}), ~~~ 
\frac{d}{d \theta^{(2)}}c(\theta_{01}, \hat{\theta}_0^{(2)}) =\bar{t}^{(2)} ,
\]
and $\hat{\theta}_0 =(\theta_{01},\hat{\theta}_0^{(2)}) $. Now, we define the signed likelihood ratio as
\[
R_t =sgn(\hat{\theta}_1 -\theta_{01}) \sqrt{2w} ((\hat{\theta} -\hat{\theta}_0 )(\bar{t}_1,\bar{t}^{(2)})^T -c(\hat{\theta}) +c(\hat{\theta}_0))^{1/2},
\]
and 
\[
U_t = \sqrt{w} (\hat{\theta}_1 -\theta_{01}) |\Sigma(\hat{\theta})|^{1/2} |\Sigma_{22}(\hat{\theta}_0)|^{-1/2},
\]
where $\Sigma(\theta) =\frac{\partial^2}{\partial\theta \partial \theta^T} c(\theta)$ and
$\Sigma_{22}(\theta) = \frac{\partial^2}{\partial\theta^{(2)} \partial \theta^{(2)T}} c(\theta)$. 
As shown in (5.2.10) of Jensen (1995), 
\[
R_t^* =R_t +\frac{1}{R_t} \ln (U_t/R_t )
\]
follows normal in the sense that $P_0(R_t^* \geq x) = (1-\Phi(x))(1+o((1+x) w^{-3/2})$. An alarm will be made at
\[
\tau_{MGL}=\inf\{ t>0: R_t^* > b \}.
\]

The evaluation of FDP will be similar to Theorem 2. We define a changed measure $P_t^{\theta} (.) $ on $X_{t-w+1},..., X_t$ such that $\theta_1 > \theta_{01} $ while $(\theta_2, ..., \theta_p)$ are unchanged (e.g. $(\theta_{02},...,\theta_{0p})$) and under $P_t^{\theta}(.)$ the log-likelihood ratio is equal to 
\[
l_t = b R_t^* -b^2 /2 ,
\]
up to the order $O(1/w^{3/2})$. At $ R_t^* $ around $b$, we use the following stochastic expansion for $R_t$. 
\begin{eqnarray*}
R_t & \approx & \sqrt{2w} ((\hat{\theta} -\hat{\theta}_0 ) (\bar{t}_1, \bar{t}^{(2)} )^T -c(\hat{\theta}) +c(\hat{\theta}_0 ))^{1/2} \\
&= & \sqrt{2w} ((\theta -\theta_0) (\bar{t}_1, \bar{t}^{(2)} )^T +(\hat{\theta}-\theta ) (\bar{t}_1, \bar{t}^{(2)} )^T -(\hat{\theta}_0-\theta_0 ) (\bar{t}_1, \bar{t}^{(2)} )^T  -c(\hat{\theta}) +c(\hat{\theta}_0 ))^{1/2} \\
& =& \sqrt{2w} ((\theta_1 -\theta_{01}) \bar{t}_1 +(\hat{\theta}-\theta ) ((\bar{t}_1, \bar{t}^{(2)} ) -\frac{d}{d\theta} c(\theta)) - (\hat{\theta}_0^{(2)} -\theta^{(2)}) (\bar{t}^{(2)} -\frac{d}{d\theta^{(2)}} c(\theta_0))^T \\
& & -\frac{1}{2} ((\hat{\theta}-\theta ) \frac{\partial^2}{\partial \theta \partial \theta^T} c(\theta) (\hat{\theta}-\theta )^T +\frac{1}{2} (\hat{\theta}_0^{(2)} -\theta^{(2)})\frac{\partial^2}{\partial \theta^{(2)} \partial \theta^{(2)T}} c(\theta_0) 
(\hat{\theta}_0^{(2)} -\theta^{(2)})^T +o_p(1/w))^{1/2}  \\
& = & \sqrt{2w} ((\theta_1 -\theta_{01})c_1'(\theta ) -c(\theta)+c(\theta_0)
+(\theta_1 -\theta_{01})(\bar{t}_1 - c_1'(\theta ) ) \\
& & +\frac{1}{2} ((\hat{\theta}-\theta ) \frac{\partial^2}{\partial \theta \partial \theta^T} c(\theta) (\hat{\theta}-\theta )^T -\frac{1}{2} (\hat{\theta}_0^{(2)} -\theta^{(2)})\frac{\partial^2}{\partial \theta^{(2)} \partial \theta^{(2)T}} c(\theta_0) 
(\hat{\theta}_0^{(2)} -\theta^{(2)})^T +o_p(1/w))^{1/2} \\
& =& \sqrt{2w \psi (\theta,\theta_0)} + (\frac{w}{2\psi (\theta,\theta_0)} )^{1/2} (\theta_1 -\theta_{01})(\bar{t}_1 - c_1'(\theta ) ) +O_p(1/w)) ,
\end{eqnarray*}
where $\psi (\theta,\theta_0) = (\theta_1 -\theta_{01})c_1'(\theta ) -c(\theta)+c(\theta_0)$ and $c_1'(\theta)$ is the first component of $dc(\theta)/d\theta$. 

The same expansion if true for $R_s$ with $s=t\pm O(\sqrt{w})$. Thus, up to the order $O(1/w)$, from $E_t^*(R_t^*) \approx b$, $\theta$ should satisfy 
\begin{equation}
\psi (\theta,\theta_0) = (\theta_1 -\theta_{01})c_1'(\theta ) -c(\theta)+c(\theta_0) = b^2/(2w) .
\end{equation}
When $b/\sqrt{w} $ is small, $\theta_1 -\theta_{01} \approx b/(w c_{11}''(\theta_0))^{1/2}$, where $c_{11}''(\theta_0) = (c''(\theta_0))_{11} $. Thus, as $ w \rightarrow \infty $, since
$l_s -l_t  \approx (\theta_1 -\theta_{01}) (\sum_{i=s-w+1}^s t_1(X_i) -\sum_{i=t-w+1}^t t_1(X_i))$, it is approximately a two-sided random walk with mean and variance
\[
E_t(l_s  -t_t) \approx bE_t(R_s -R_t ) \approx - (\theta_1 -\theta_{01}) |t-s| (c_1'(\theta) -c_1'(\theta_0 )),
\]
\[
Var_t(l_s -l_t) \approx ( \theta_1 -\theta_{01})^2 |t-s| (c_{11}''(\theta_0 )+c_{11}''(\theta)) .
\]
Thus, the FDP can be approximated
\[
P_0^*(\tau_{MLR} \leq L ) \approx \frac{L(\theta_1-\theta_{01}) (c_1'(\theta)-c_1'(\theta_0))}{b} \phi(b) \upsilon ((\theta_1-\theta_{01})  \sqrt{c_{11}''(\theta_0)+c_{11}''(\theta)} )  \]
\begin{equation}
\approx \frac{Lb}{w} \phi(b ) e^{-b\rho_+/\sqrt{w c_{11}''(\theta_0)}},
\end{equation}
where $\rho_+$ will be obtained under $P_{\theta_0}^*(.)$ by using the random walk $S_t = \sum_{k=1} (t_1(X_k)- t_1(X_k')$ with $X_k'$ being independent copy of $X_k$. \qed

\noindent{\bf Remark:} Equation (4.8) may have multiple solutions, we can take $\theta^{(2)}$ to be any specific convenient values. Note that in multiparameter case, $U_t$ is different from Wald' statistic and Rao's score statistic. Wald's statistic is
\[
W_t = \sqrt{w} (\hat{\theta}_1 -\theta_{01}) |\sigma_{11} -\Sigma_{12} \Sigma_{22}^{-1} \Sigma_{21}|^{1/2} ,
\]
where all $\Sigma$'s are evaluated at $\hat{\theta}$. Rao's statistic, however, is defined as
\[
V_t = \sqrt{w}  (\bar{t}_1 -c_{\theta_1}'(\hat{\theta}_0))  |\sigma_{11} -\Sigma_{12} \Sigma_{22}^{-1} \Sigma_{21}|^{-1/2},
\]
where all $\Sigma$'s are evaluated at $\hat{\theta}_0$. 

Similarly, we can generalize the CUSUM and S-R procedures to the windowed version by estimating the nuisance parameter $\theta^{(2)}$ as $\hat{\theta}_0^{(2)}$ based on the data in the window at $\theta_{01}$. Specifically, suppose we want to detection an increase from $\theta_{01}$ to a reference value $\theta_{11}$. Denote by $\hat{\theta}_k^{(i)} =(\theta_{i1}, \hat{\theta}_i^{(2)})$ as the corresponding MLE for $\theta$ based on the data $X_{t-k+1}, ..., X_t$, for $i=0,1$. For $k=1, ..., w$, we can define the profile log-likelihood ratio as
\[
l_k(t) = k[(\hat{\theta}_k^{(1)} -\hat{\theta}_k^{(0)}) \bar{t}_{t;k} -c(\hat{\theta}_k^{(1)}) +c(\hat{\theta}_k^{(0)})].
\]
Then the CUSUM procedure is defined to make an alarm at
\begin{equation}
\tau_{CS} =\inf\{ t>0: \max_{1\leq k\leq w} l_k(t) >d \}.
\end{equation}
Similarly, the S-R procedure is defined to make an alarm at
\begin{equation}
\tau_{SR} =\inf \{ t>0: \ln \sum_{1\leq k\leq w} \exp(l_k(t)) >c\} .
\end{equation}
For the above window restricted CUSUM and S-R procedures, it is often better to replace $\hat{\theta}_k^{(i)}$ by $\hat{\theta}_w^{(i)}$ as shown in the following example.

\noindent{\bf Example 3.} We consider the detection of variance $\sigma^2$ increase from 1 when the mean $\mu$ is treated unknown in normal random variable. By writing $ \theta_1 = -\frac{1}{2\sigma^2}$, $\theta_2 =\mu /\sigma^2$, $t_1(x) = x^2$, and $ t_2(x) =x$, the density function of $N(\mu, \sigma^2)$ has an exponential form with
\[
c(\theta_1, \theta_2) = -\frac{\theta_2^2}{4\theta_1}- \frac{1}{2} \ln(-2\theta_1) . 
\]
For observations $X_{t-w+1},..., X_t$, we denote by $\bar{X}_{t;w} = \sum_{i=1}^w X_{t-i+1}/w$ and 
$\sigma_{t;w}^2 =\sum_{i=1}^w (X_{t-i+1}- \bar{X}_{t;w})^2 /w$. 

From
\[
\frac{d}{d\theta}c(\theta) =(\frac{\theta_2^2}{4\theta_1^2}-\frac{1}{2\theta_1}, -\frac{\theta_2}{2\theta_1}) =(\bar{t}_1, \bar{t}_2) =(\sum_{i=1}^w X_{t-i+1}^2/w, \bar{X}_{t;w}),
\]
we get
\[
\hat{\theta}_1 = -\frac{1}{2\hat{\sigma}_{t;w}^2}; ~~~~ \hat{\theta}_2 = \frac{\bar{X}_{t;w}}{\hat{\sigma}_{t;w}^2}.
\]
This gives $l(\hat{\theta}) = (\hat{\theta}_1, \hat{\theta}_2) ( \bar{t}_1, \bar{t}_2)^T -c(\hat{\theta}_1, \hat{\theta}_2) = -w/2 -w\ln (\hat{\sigma}_{t;w}^2) /2 $. 

At $\theta_{01} =-1/2 (\sigma^2=1 )$, $\hat{\theta}_{02} = \bar{X}_{t;w}$. This gives $l(\theta_{01}, \hat{\theta}_{02} ) = -w \hat{\sigma}_{t;w}^2 /2 $. 

Thus, 
\[
R_t = sgn(\hat{\sigma}_{t;w}^2 -1) \sqrt{w} (\hat{\sigma}_{t;w}^2 -1-\ln \hat{\sigma}_{t;w}^2 )^{1/2}.
\]

On the other hand, we can find 
\[
|\Sigma(\theta)| =|\frac{\partial^2}{\partial\theta \partial \theta^T} c(\theta)| = 1/(4\theta_1^3 ),
\]
and $|\Sigma_{22}(\theta)|  = -1/(2\theta_1 )$.  Thus
\[
U_t =\sqrt{w/2} (\hat{\sigma}_{t;w}^2 -1) \hat{\sigma}_{t;w}.
\]
Finally, we have
\[
R_t^* = R_t +\frac{1}{2R_t} \ln \frac{\hat{\sigma}_{t;w}^2 (\hat{\sigma}_{t;w}^2 -1)^2/2}{\hat{\sigma}_{t;w}^2 -1 -\ln \hat{\sigma}_{t;w}^2} .
\]
Note that in this example, $W_t =V_t =\sqrt{w/2} (\hat{\sigma}_{t;w}^2 -1) $. 

On the other hand, for given $\theta_1 =\theta_{i1}$ ( $\sigma^2 =\sigma_i^2$) and $X_{t-k+1},..., X_t$, $\hat{\theta}_2^{(i)} =-2\theta_{i1} \bar{X}_{t;k}$ for $i=0,1$. A simplification shows that
\[
l_k(t) = \frac{k}{2} [(\frac{1}{\sigma_0^2} -\frac{1}{\sigma_1^2}) (\hat{\sigma}_{t;k}^2 +\ln (\sigma_0^2/\sigma_1^2) ].
\]
The CUSUM and S-R procedures can be defined based on (4.6) and (4.7). In this example, since the nuisance parameter is the mean for given $\sigma^2$, a better estimate for $\mu$ is $\bar{X}_{t;w}$ when the data from the whole window are available. In fact, simulation results show that it gives slightly larger POD. 

Table 3 gives the comparison of the five charts based on $W_t$, $R_t$, $R_t^*$, CUSUM and S-R charts similar to Table 2 with same parameter. For the CUSUM and S-R procedures, the reference value for the variance after change is taken as 2.   The numbers in the brackets are the approximations for the FDP. The chart based on $R_t^*$ gives much better results. 
Again, we see that $W_t$ has much heavier tail and $R_t$, somehow, has lighter tail that $R_t^*$.

\begin{table}[ht]
\begin{center}
Table 3.  Comparison of POD for $L=20$ and $w=20$ for variance increase in normal R.V. with unknown mean \\
\begin{tabular}{|c |c |c |c |c| c|}   \hline
$\sigma^2$ & $W_t$ & $R_t$ & $R_t^*$ &CUSUM &S-R \\
boundary & b=3.05 & b=2.40 & b=2.65 &d=3.58& c= 5.45\\ \hline
1.00 & 0.0202(0.0066) & 0.0218(0.0345) & 0.0207(0.0194)&0.0204 &0.0206 \\
1.25 & 0.0722 & 0.0731 & 0.0721&0.0798 &0.0795 \\
1.50 & 0.1833 & 0.1861 & 0.1810 &0.1939 &0.1986\\
1.75 & 0.3327 & 0.3433 & 0.3287&0.3430 & 0.3462\\
2.00 & 0.4797 & 0.4819 & 0.4808 &0.4937& 0.5064\\
2.25 & 0.6118 & 0.6192 & 0.6183 &0.6264& 0.6333\\
2.50 & 0.7335 & 0.7311 & 0.7259 &0.7333& 0.7350\\
2.75 & 0.8157 & 0.8100 & 0.8101 &0.8138& 0.8149 \\
3.0 & 0.8659 & 0.8673 & 0.8601 &0.8671 & 0.8688 \\ \hline
\end{tabular}
\end{center}
\end{table}

\subsection{Detection of a scale parameter change as a curved exponential family}
A discussion for the likelihood ratio based test in the two-parameter exponential family case is given in Siegmund and Yakir (2000), where a mixture is taken for the nuisance parameter. 

Here, we consider the problem as detecting the change of an interest scale parameter, not necessarily a component of the canonical parameters. For example, in the normal model $N(\mu, \sigma^2 )$, we may be interested in detecting the change of $\mu$ or $\sigma$ while the other parameter may be subject to change, and neither of the two is a canonical parameter. 

 For this purpose, we assume the canonical parameter $\theta =(\theta_1, ..., \theta_p)$ is a one-to-one differentiable  transformation, written as $\theta =\theta (\beta)$ from the original parameter $\beta = (\beta_1, ...,\beta_p) $. 
 We are interested in detecting a change of $\beta_p $ from $\beta_{0p}$, while $\beta^{p-1} =(\beta_1 ..., \beta_{p-1})$ are assumed unknown, and  may be randomly or slowly varying. 
 
 Denote $\hat{\beta}$ as the MLE by solving 
 \[
 (\bar{t} -c'(\theta)) \frac{d\theta}{d\beta^T} (\beta) =0,
 \]
 and  for $\beta_0 =(\beta_1, \beta_2, ..., \beta_{0p})$, $\hat{\beta}_0$ be solution of 
 \[
(\bar{t} -c'(\theta)) \frac{d\theta}{d\beta^{p-1T}} (\beta_0) =0. 
\]
By writing $\hat{\theta} =\theta(\hat{\beta})$ and $\hat{\theta}_0 =\theta(\hat{\beta}_0)$, the signed likelihood ratio can be written as
\[
R_t =\sqrt{2w} sgn(\hat{\beta}_p -\beta_{0p}) ((\hat{\theta} -\hat{\theta}_0) \bar{t} -c(\hat{\theta}) +c(\hat{\theta}_0))^{1/2}.
\]
The form of $U_t$ is quite complicated and we use a form that is given as (5.3.19) of Jensen (1995, Section 5.3):
\[
U_t=\sqrt{w} sgn(R_t) |\Sigma(\hat{\theta})|^{1/2} |j^{p-1}(\hat{\beta}^{p-1})|^{-1/2} |J| ,
\]
where 
\[
J=(\frac{\partial\theta^T}{\partial\beta^{p-1}} (\hat{\beta}_0), (\hat{\theta} -\hat{\theta}_0)^T), 
\]
is the Jacobian and
\[
|j^{p-1}(\hat{\beta}^{p-1})| =\prod_{i=1}^{p-1} |j_{22}^i - j_{21}^i(j_{11}^i)^{-1}j_{12}^i|,
\]
and $j_{kl}^i = (-\frac{\partial^2 l}{\partial \beta^i \partial \beta^{iT}} (\hat{\beta}^{i-1} ))_{kl}$ with $l(\beta )= \theta \bar{t} -c( \theta(\beta))$.

When $\beta =\theta$ as the canonical parameter, $J =(I_{p \times (p-1)}, (\hat{\theta} -\hat{\theta}_0 )^T)$ and $|J| = |\hat{\theta}_p -\hat{\theta}_{0p}| $ as discussed before.

The adjusted signed likelihood ratio test statistic $R_t^* =R_t +\ln(U_t/R_t)/R_t$ can then be treated as a second order corrected normal random variable. An alarm will be made at
\[
\tau =\inf\{t>w: R_t^* >b \}
\]

To derive an approximation for FDP, we follow the same lines as in the canonical parameter case. Define a changed measure $P_t^{\beta}(.)$ with $\theta =\theta (\beta)$ on $X_{t-w+1},..., X_t$ such that  under $P_t^{\theta}(.)$ the log-likelihood ratio is equal to 
\[
l_t = b R_t^* -b^2 /2 ,
\]
up to the order $O(1/w^{3/2})$. At $ R_t^* $ around $b$, the following stochastic expansion for $R_t$ can be derived as before
\[
R_t \approx \sqrt{2w \psi (\theta,\theta_0)} + (\frac{w}{2\psi (\theta,\theta_0)} )^{1/2} (\theta -\theta_0)(\bar{t} - \frac{d}{d\theta^T}c(\theta ) ) +O_p(1/w)) ,
\]
where $\psi (\theta,\theta_0) = (\theta-\theta_0)c'(\theta ) -c(\theta)+c(\theta_0)$ with $c'(\theta) = \frac{d}{d\theta^T}c(\theta )$ and $c''(\theta)= \frac{\partial^2}{\partial \theta \partial \theta^T} c(\theta)$. 

The same expansion is true for $s=t\pm O(\sqrt{w})$. Thus, up to the order $O(1/w)$, from $E_t^*(R_t^*) \approx b$, $\theta$ should satisfy 
\[
\psi (\theta,\theta_0) = (\theta-\theta_0)c'(\theta ) -c(\theta)+c(\theta_0) = b^2/(2w)
\]
Similarly, $l_s-l_t \approx (\theta-\theta_0)(\sum_{i=s-w+1}^s t(X_i) -\sum_{i=t-w+1}^t t(X_i) )$ can be approximated by a random walk with drift and variance
\[
E_t[l_s-l_t] \approx - (t-s) (\theta -\theta_0)(c'(\theta)-c'(\theta_0)),
\]
\[
Var_T(l_s -l_t) \approx  (t-s) (\theta -\theta_0)((c''(\theta)+c''(\theta_0))(\theta -\theta_0)^T.
\]

The FDP can be approximated by 
\[
P_0^*(\tau \leq L) \approx \frac{L (\theta -\theta_0)(c'(\theta)-c'(\theta_0)) }{ b} \phi(b) \upsilon ( ((\theta -\theta_0)((c''(\theta)+c''(\theta_0))(\theta -\theta_0)^T)^{1/2}) 
\]
\[
\approx \frac{Lb}{w} \phi(b) e^{-\rho_+ b/\sqrt{w}} ,
\]
where $\rho_+$ will be calculated as $\rho_+ =E_{\theta_0} S_{\tau_+}^2/(2E_{\theta_0}S_{\tau_+})$ from the ladder time defined by
\[
\tau_+ =\inf\{ t>0: S_t =\gamma (\sum_{i=1}^t (t(X_i)-t(X_i'))>0 \},
\]
where $X_i$ follows $f_{\theta_0}(.)$  and $X_i'$ is independent copy of $X_i$ and 
$\gamma =\lim_{\theta \rightarrow \theta_0} (\theta-\theta_0)/((\theta -\theta_0)c''(\theta_0)(\theta -\theta_0)^T)^{1/2}$ is the direction. \qed

In principle, the above technique can be applied to any multi-parameter parametric family. However, the calculation of $U_t$ in general setting is generally complicated beyond the exponential family, except in some special cases, e.g. when $\beta_p$ is orthogonal to $\beta^{p-1}$ (see  Barndorff-Nielsen and Cox (1994, Section 6.6) or  Severini (2000, Section 7.5)). A different way of deriving $U_t$, developed by Don Fraser, is explained in Section 4.3 of Davison and Reid (2022), where a tangent exponential model is also discussed for general parametric family. A recent discussion on the approximation for the significance function with increasing number of nuisance parameters is given in Tang and Reid (2020). 

However, the CUSUM and S-R procedures can be used like the one-parameter case.  Let $f(x|\beta)$ be the density function and denote
\[
l_k^t (\beta ) = \sum_{i=t-k+1}^t f(X_i|\beta).
\]
To detect a change of $\beta_p$ from $\beta_{0p}$ to a reference value $\beta_{1p}$, let
$\hat{\beta}_k^{(i)}$ be the MLE of $\beta$ given $\beta_p =\beta_{ip}$ for $i=0,1$. Then the CUSUM and S-R charts can be defined to make an alarm at
\[
\tau_{CS} =\inf\{t>0: \max_{1\leq k\leq w} (l_k^t (\hat{\beta}_k^{(1)} )-l_k^t (\hat{\beta}_k^{(0)})) > d\},
\]
\[
\tau_{SR} = \inf\{t>0: \ln(\sum_{1\leq k\leq w}\exp(l_k^t (\hat{\beta}_k^{(1)})-l_k^t( \hat{\beta}_k^{(0)}))) > c\}.  
\]
When the window restricted data are available, as shown in the following example, it is often better to replace $\hat{\beta}_k^{(i)}$ by $\hat{\beta}_w^{(i)}$. 

\noindent{\bf Example 4.} We consider Example 3 with $ \theta_1 = -\frac{1}{2\sigma^2}$, $\theta_2 =\mu /\sigma^2$, $t_1(x) = x^2$, and $ t_2(x) =x$, the density function of $N(\mu, \sigma^2)$ has an exponential form with
\[
c(\theta_1, \theta_2) = -\frac{\theta_2^2}{4\theta_1}- \frac{1}{2} \ln(-2\theta_1) . 
\]
Here we detect change of $\mu$ from 0 to a positive number $\delta$ with $\beta =(\sigma^2, \mu)$. Again, $l(\hat{\theta}) = (\hat{\theta}_1, \hat{\theta}_2) ( \bar{t}_1, \bar{t}_2)^T -c(\hat{\theta}_1, \hat{\theta}_2) = -w/2 -w\ln (\hat{\sigma}_{t;w}^2) /2 $. At $\mu =0$, $\hat{\theta}_{01} = -1/(2\hat{\sigma}_0^2)$ for
$\hat{\sigma}_0^2 =\bar{t}_1 =\sum_{i=t-w+1}^t X_i^2 /w$ and 
\[
l(\hat{\theta}_0) = -w/2-w\ln(\hat{\sigma}_0^2) /2.
\]
Thus, 
\[
R_t =sgn(\bar{X}_{t;w}) \sqrt{w} (\ln ( \hat{\sigma}_0^2/\hat{\sigma}_{t;w}^2 ))^{1/2}.
\]
To calculate $U_t$, we first note that $ |\Sigma(\hat{\theta})| = 1/(4\hat{\theta_2}^3)$. This gives $|\Sigma(\hat{\theta})|^{1/2} = \sqrt{2}\hat{\sigma}_{t;w}^3  $.

Second, we can calculate 
\[
|J| = \bar{X}_{t;w}/ (2 \hat{\sigma}_0^4 \hat{\sigma}_{t;w})^2 .
\]
Third, since $p=2$, 
\[
|j^{p-1}(\hat{\beta}^{p-1})| =\frac{\partial^2 l}{(\partial \sigma^2)^2} (\hat{\beta}_0) = 1/(2\hat{\sigma}_0^4).
\]
A simplification gives
\[
U_t =\sqrt{w}\bar{X}_{t;w} \hat{\sigma}_{t:w}/\hat{\sigma}_0^2 .
\]
Thus,
\[
R_t^* =R_t +\frac{1}{R_t} \ln (\frac{|\bar{X}_{t;w}| \hat{\sigma}_{t:w}/\hat{\sigma}_0^2}{\ln ( \hat{\sigma}_0^2/\hat{\sigma}_{t;w}^2 ))^{1/2}}).
\]

For the CUSUM and S-R procedure, for the reference $\mu_1 =\delta$, we use $\hat{\sigma}_{t;k}^{(i)2} =\sum_{j=t-k+1}^t (X_j -\mu_i)^2 /k $ as the estimator for $\sigma^2$ for $i=0,1$. Thus, the CUSUM procedure is defined as to make an alarm at
\[
\tau_{CS}=\inf\{t>0: \max_{1\leq k \leq w}  \frac{k}{2} \ln \frac{\hat{\sigma}_{t;k}^{(0)2}}{\hat{\sigma}_{t;k}^{(1)2}} >d\}
\]
and the S-R procedure makes an alarm at
\[
\tau_{SR}=\inf\{ t>0: \ln \sum_{1\leq k \leq w} \left(\frac{\hat{\sigma}_{t;k}^{(0)}}{\hat{\sigma}_{t;k}^{(1)}}\right)^k > c\}
\]
However, the simulation results show that the POD is rather poor. Instead, we use 
$\hat{\sigma}_{t;w}^2$ as a pooled estimate for $\sigma$ that is free of $\mu$. The CUSUM chart changes to make an alarm at
\[
\tau_{CS}=\inf\{t>0: \max_{1\leq k\leq w} \sum_{t-k+1}^t (\delta X_j -\delta^2/2)/\hat{\sigma}_{t;w}^2 > d \}
\]
and the S-R procedure becomes to make an alarm at
\[
\tau_{SR} =\inf\{ t>0: \ln(\sum_{1\leq k\leq w} \exp(\sum_{t-k+1}^t (\delta X_j -\delta^2/2)/\hat{\sigma}_{t;w}^2)) > c\}.
\]

Table 4 gives a comparison between $R_t$, $R_t^*$, the regular $t$-test based on $\sqrt{w} \bar{X}_{t:w} /\hat{\sigma}_{t;w}$, CUSUM and S-R charts. For the CUSUM and S-R procedure the reference value is taken as $\delta =0.5$.  We see that both the $t-$test and $R_t$ test have heavier tails than $R_t^*$ test.  Also, in this example, the CUSUM and S-R procedure do not perform as well as $R_t^*$.

\begin{table}[ht]
\begin{center}
Table 4.  Comparison of POD for $L=20$ and $w=20$ for mean increase in normal R.V. with unknown variance \\
\begin{tabular}{|c |c |c |c |c|c|}   \hline
$\mu$ & $t$-test & $R_t$ & $R_t^*$ &CUSUM &S-R\\
boundary & b=3.10 & b=2.77 & b=2.67 & d=4.26 &c=6.7\\ \hline
0.00 & 0.0234(0.0066) & 0.0246(0.0143) & 0.0237(0.0194)& 0.0233 &0.0240 \\
0.25 & 0.0969 & 0.1099 & 0.1088&0.0965& 0.1033 \\
0.5 & 0.3546 & 0.3576 & 0.3647 & 0.3156 &0.3361 \\
0.75 & 0.7103 & 0.7159 & 0.7210 & 0.6400 &0.6686 \\
1.00 & 0.9338 & 0.9402 & 0.9378 & 0.8816 & 0.8968\\
1.25 & 0.9933 & 0.9931 & 0.9939 & 0.9766 & 0.9798  \\
1.50 & 0.9999 & 0.9996 & 0.9997 & 0.9973 & 0.9978\\
1.75 & 1.0 & 1.0 & 1.0 &0.9998 & 0.9998 \\
2.0 & 1.0 & 1.0 &  1.0 &1.00 & 1.0000\\ \hline
\end{tabular}
\end{center}
\end{table}

\noindent{\bf Example 5.}  We use IBM stock daily closing prices from May 2, 2021 to May 2, 2022 (downloaded from finance.yahoo.com) to demonstrate the applications. The simple geometric random walk model is used. By taking the differences of the logarithms, the first figure in Figure 1 shows the plot the data. To make the data stable, we first trimmed the outliers using the $\pm 3 \sigma$-rule. The auto-correlation functions of the data and squared data do not show significant correlations. 

\begin{figure}
\includegraphics[width =\textwidth,height=9.5in]{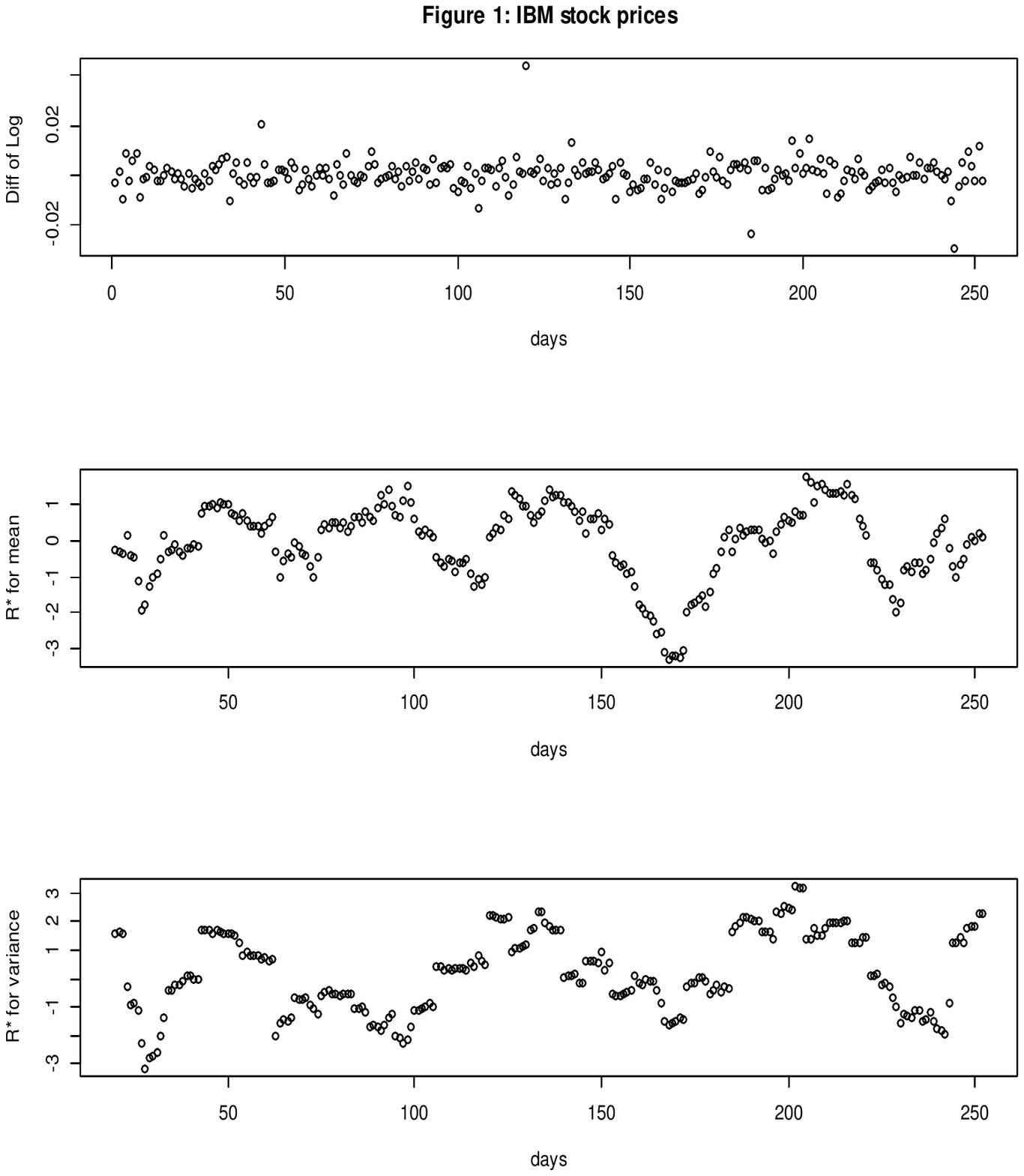}
\end{figure}

We take $w=20$. So the lag-time is 20. The second figure gives the plot of $R_t^*$ for detecting the change of mean based on Example 4. No alarm is raised for positive change of mean. However, we note that the chart can also be used to detect negative change mean (b=-2.67). We found an alarm is raised at day 167 with $R^* =-3.094$ and continued until day 172. 

For detecting the variance increase, we used the standard deviation from the first 100 days data to standardize the data. the third figure gives the plot of $R^*$. An alarm is raised at day 202 (b=2.65) continued for 3 days. (We also note that there is also an alarm for variance decrease at around day 30.)

\section{Detecting Change in Multi-parameter Case}

In this communication, we considered the detection of one parameter change in one or multi-parameter exponential family. We propose to use the adjusted signed likelihood ratio as the detecting rule by using its second order normality. The advantage is that the alarm limit is quite stable and easy to design. The method introduced here can also be applied to detection of two-sided change or generally, multi-parameter change. For example, we can combine two charts for detecting possible change in either mean or variance under the normal model. Generally, without specifying the direction of change, the Bartlett adjustment can be used to find a higher order chi-square approximation.  
We briefly give the main steps and leave the details for further discussions. 

For $p$-dimensional parameter $\theta$, denote $l_{t;w} (\theta; \hat{\theta}, a) =\sum_{k=t-w+1}^t \ln f_{\theta}(X_k) $ as the log-likelihood where $(\hat{\theta}, a)$ is a transform of the data and $a$ is asymptotic ancillary, and $\hat{\theta}_t$ is the MLE.  For detecting change from $\theta =\theta_0$, the log-likelihood ratio statistic can be written as 
\[
W_t =2(l_{t;w}(\hat{\theta}) - l_{t;w}(\theta_0)). 
\]
As $w\rightarrow \infty$, let
\[
E_{\theta_0} W_t = p(1+b(\theta_0)/w +O(1/w^2)).
\]
Then the Bartlett adjusted $W_t$ defined as
\[
W_t^* = W_t/(1+b(\theta_0)/w),
\]
follows chi-square distribution up to the order $O(w^{-2})$. 
Systematic discussions on the multi-parameter case can be seen in Barndorff-Nielsen and Cox (1984) and explicit formulas for $b(\theta_0)$ are given in Section 6.5.1 of Barndorff-Nielsen and Cox (1994) or Section 8.8 of Brazzale, et al. (2006). 
 An alarm will be raised at
\[
\tau_W = \inf\{ t>0: W_t^* >b^2 \} 
\]
The following theorem gives the approximation for the FDP. 

\noindent{\bf Theorem 3.} {\it As $b, w \rightarrow \infty$ such that $b/\sqrt{w} \rightarrow h>0$, 
\begin{equation}
P_0^*(\tau_W \leq L) \approx  L \frac{(b^2/p)^{p/2} e^{p/2}}{2\sqrt{p\pi}} e^{-b^2/2} 
(\theta -\theta_0)(c'(\theta)-c'(\theta_0)) \upsilon(\theta^* ([(\theta-\theta_0) (c''(\theta) + c''(\theta_0))(\theta -\theta_0)^T]^{1/2} ). 
\end{equation}
}

{\em Proof.} The evaluation of FDP will be carried by directly assuming $W_t^*$ follows chi-square distribution with $p$ degrees of freedom. Define a changed measure $P_t^*(.)$ on $X_{t-w+1},..., X_t$ such that under $P_t^*(.)$ the log-likelihood ratio of $X_1,..., X_t$ under $P_t^*(.)$ with respect to $P_0^*(.)$ is equal to
\[
l_t =\theta^* W_t^* +\frac{p}{2} \ln (1-2\theta^*).
\]
The value of $\theta^*$ is taken such that
\[
E_t^*[W_t^*] \approx \frac{d}{d\theta^*} [-\frac{p}{2} \ln (1-2\theta^*)] = \frac{p}{1-2\theta^*} = b^2.
\]
This gives $\theta^*= \frac{1}{2} (1-\frac{p}{b^2})$. Now let $\tilde{l}_t =\theta^* (W_t^* -b^2)$. So $\tilde{l}_t$ has approximately mean zero and variance $2p \theta^{*2}/(1-2\theta^*)^2$. 

By assuming $w\rightarrow \infty$ and $b\rightarrow \infty $, we can have the following approximation for FDP 
\begin{eqnarray*}
P_0^*(\tau_W \leq L) &=& P_0^*(\max_{1\leq t\leq L} W_t^* >b^2) \\
&\approx& L e^{-p(\theta^*/(1-2\theta^*) +\ln (1-2\theta^*)/2)} \frac{1}{\sqrt{2\pi}} \frac{1-2\theta^*}{\sqrt{2p}\theta^* } \lim_{L\rightarrow \infty}E_t^* [\frac{M_t}{S_t}] \\
&\approx & L \frac{(b^2/p)^{p/2} e^{p/2}}{\sqrt{p\pi}(1-p/b^2)} e^{-b^2/2} \lim_{L\rightarrow \infty}E_t^* [\frac{M_t}{S_t}] \\
&\approx & L \frac{(b^2/2)^{p/2-1} e^{-b^2/2}}{\Gamma(p/2)(1-p/b^2)} \lim_{L \rightarrow \infty}E_t^* [\frac{M_t}{S_t}] ,
\end{eqnarray*}
where the last equation is obtained by using Sterling's formula when $p$ is relatively large and $M_t =\max_{1\leq s\leq L} (l_s-l_t)$ and $S_t =\sum_{1\leq s \leq L} (l_s-l_t) $. Similar approximation is also obtained in Wu and Siegmund (2022) for the multi-dimensional data stream case. 

In the following, we just need to show that $l_s-l_t$ is approximately a local random walk. For simplicity, we only verify it in the exponential family case by assuming that $X_{t-w+1}, ...., X_t$ follow $f_{\theta}(x) =\exp(\theta t(x) -c(\theta))$. By ignoring the lower order terms, we have
\begin{eqnarray*}
W_t^* \approx W_t &=& 2w [(\hat{\theta}_t -\theta_0 ) \bar{t} -c(\hat{\theta}_t)+c(\theta_0)] \\
&\approx &2w[(\theta -\theta_0) \bar{t} -c(\theta)+c(\theta_0) +\frac{1}{2} (\hat{\theta}_t -\theta)c''(\theta) (\hat{\theta}_t -\theta)^T ].
\end{eqnarray*}
Similar expansion is true for $s=t \pm O(\sqrt{w}) $ under $P_t^*(.)$. Thus, at the first order, $\theta$ should be chosen such that
\[
2w((\theta-\theta_0)c'(\theta)-c(\theta)+c(\theta_0) = b^2.
\]
Now it is easy to see that $l_s-l_t = \theta^*(W_s^*-W_t^*)$ behaves like a two sided random walk with drift $-\theta^*|t-s| (\theta -\theta_0)(c'(\theta)-c'(\theta_0)) $ and variance $ \theta^{*2} |t-s| (\theta-\theta_0) ( c''(\theta) +c''(\theta_0))(\theta -\theta_0)^T $.
Therefore
\[
\lim_{L \rightarrow \infty}E_t^* [\frac{M_t}{S_t}] \approx 
\theta^* (\theta -\theta_0)(c'(\theta)-c'(\theta_0)) \upsilon(\theta^* ([(\theta-\theta_0) (c''(\theta) + c''(\theta_0))(\theta -\theta_0)^T]^{1/2} ). 
\]
The result is obtained by combining the above results. When $b^2/w$ is small, $\theta -\theta_0$ is also small. The approximation can be simplified to 
\begin{eqnarray}
P_0^*(\tau_W \leq L) &\approx& \frac{Lb^2 (b^2/p)^{p/2}e^{p/2}}{4w \sqrt{p\pi}} e^{-b^2/2}  \upsilon ((1-p/b^2)(\theta-\theta_0) (c''(\theta) +c''(\theta_0))(\theta -\theta_0)^T]^{1/2}/\sqrt{2} ) \nonumber \\
&\approx& \frac{Lb^2 (b^2/p)^{p/2}e^{p/2}}{4w \sqrt{p\pi}} e^{-b^2/2} e^{-\rho_+ b(1-p/b^2)/\sqrt{2w}}, 
\end{eqnarray}
where $\rho_+$ is calculated as $\rho_+ =E_{\theta_0}S_{\tau_+}^2/(2E_{\theta_0}S_{\tau_+})$ from the ladder time defined by
\[
\tau_+ =\inf\{ t>0: S_t= \gamma (\sum_{i=1}^t (t(X_i)-t(X_i'))>0 \},
\]
where $X_i$ follows $f_{\theta_0}(.)$  and $X_i'$ is independent copy of $X_i$ and 
$\gamma =\lim_{\theta \rightarrow \theta_0} (\theta-\theta_0)/((\theta -\theta_0)c''(\theta_0)(\theta -\theta_0)^T)^{1/2}$ is the direction. \qed

We use the normal model in Example 3 as illustration.

\noindent{\bf Example 6.} We consider the detection of change of both mean $\mu$ and variance $\sigma^2$ from 0 and 1 respectively. The likelihood ratio statistics for $H_0:\mu=0 , \sigma^2=1$ vs $H_a: \mu \neq 0$ or $\sigma^2 \neq 1 $ is 
\[
W_t = 2(l(\hat{\theta})-l(\theta_0)) = w (\bar{X}_{t;w}^2 +\hat{\sigma}^2 -1 -\ln \hat{\sigma}^2 ),
\]
where $\hat{\sigma}^2 = \sum_{t-w+1}^t (X_i -\bar{X}_{t;w})^2 /w $. At the first order, $W_t$ is asymptotically chi-square with 2 degrees of freedom. It is easy to calculate 
\[
E_0[W] =2 (1+ \frac{3}{4w} +o(1/w)).
\]
The Bartlett's adjustment gives
\[
W_t^* = W_t/(1+3/(4w)) ,
\]
with error at the order of $O(w^{-2})$. An alarm is made at the stopping time
\[
\tau_W =\inf\{ t>0: W_t^* >b^2 \}
\]
Table 5 gives simulation results for $\mu $=0, 0.25, 0.5, 0.75, 1.0, 1.25, 1.5, 1.75, 2.0 
and $\sigma^2$ =1, 1.25, 1.5, 1.75, 2.0, 2.25, 2.5, 2.75, 3.0. To match with Tables 3 and 4, we take FDP=0.08 covering four one-sided detection charts. The control limit is found to be $b^2 \approx 9.3$. As a numerical comparison, the approximation without the overshoot correct in (5.13) gives FDP = 0.112; while by multiplying the overshoot correct with normal approximation gives FDP=0.0821.  The simulation is replicated 50,000 times. 

\begin{table}[ht]
\begin{center}
Table 5.  POD for $L=20$ and $w=20$ for detecting both mean and variance change in normal R.V. \\
\begin{tabular}{|c |c |c |c |c|c|c|c|c|c|}   \hline
$\mu$/$\sigma^2$ & 1 &1.25 &1.5 &1.75 &2.0 &2.25& 2.5 &2.75&3.0\\  \hline
0.00 & 0.0807 &0.1926 &0.3270&0.4638&0.5862 &0.6873&0.7699&0.8261&0.8746 \\
0.25 & 0.0986 &0.2165 &0.3514&0.4893&0.6087 &0.7066&0.7796&0.8414& 0.8816\\
0.5 & 0.1661 & 0.3041 &0.4456&0.5735&0.6784 &0.7619&0.8239&0.8700&0.9076 \\
0.75 & 0.2833 & 0.4405&0.5802&0.6870&0.7720 & 0.8319&0.8786&0.9132&0.9355 \\
1.00 & 0.4617 &0.6141 &0.7221&0.8104&0.8645 &0.9030&0.9291&0.9487&0.9638 \\
1.25 & 0.6519 &0.7782 &0.8530&0.9029&0.9323 &0.9537&0.9661&0.9762& 0.9823 \\
1.50 & 0.8216 &0.8982 &0.9364&0.9591&0.9736 &0.9823&0.9867&0.9908&0.9932 \\
1.75 & 0.9363 &0.9660 &0.9805&0.9876& 0.9913 &0.9947&0.9961&0.9970&0.9979 \\
2.0 & 0.9847&0.9926 & 0.9958& 0.9973&0.9982 & 0.9988&0.9987&0.9992& 0.9993 \\ \hline
\end{tabular}
\end{center}
\end{table}

\noindent{\bf Example 5.} (Cont'd) For the IBM stock prices considered in Example 5, Figure 2 gives the plot of the  adjusted likelihood ratio statistic with window width 20. We can see in addition to the three changed segments detected from the four one-sided charts, another change is detected around day 100 that is caused by combination of mean increase and variance decrease.

\begin{figure}[ht]
\includegraphics[width =\textwidth,height=5in]{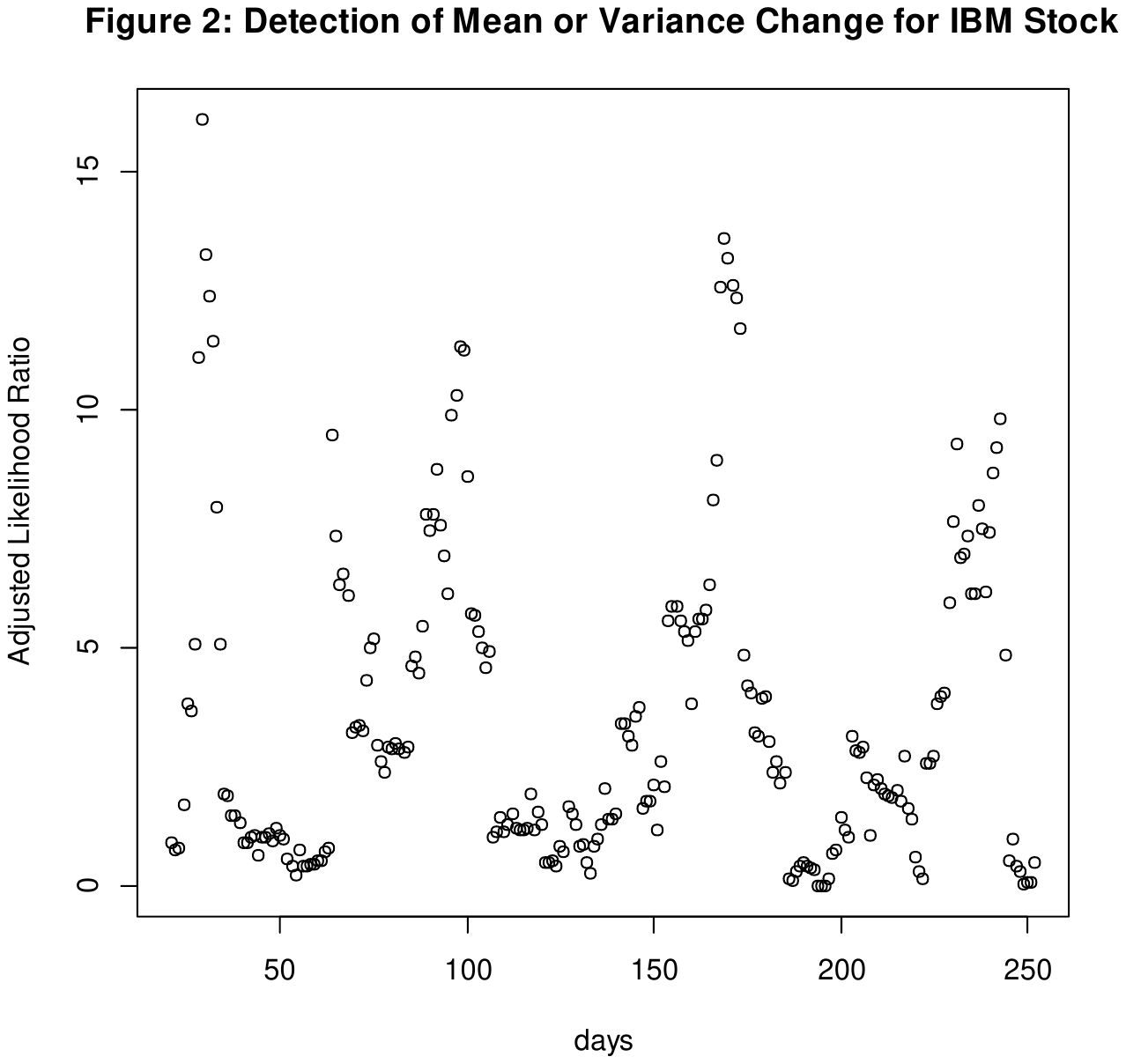}
\end{figure}

Finally we should mention that the detection of common transient signal in multidimensional data streams will be treated as a different topic as cross-sectional dependence appears. The detection of transient mean change in multivariate normal data is discussed in Wu and Siegmund (2022). The details for more general models will be left as future communications.

\section*{Acknowledgement} This was done when the author was visiting Department of Statistics of Stanford University. Many discussions with Professor D. Siegmund and comments from Professor N. Reid are greatly appreciated. 

\section*{ References}

\begin{description}
\item[{\rm Barndorff-Nielsen, O.E., Cox, D. R., 1984.}] Bartlett adjustments to the likelihood ratio statistic and the distribution of the maximum likelihood estimator. {\em Journal of Royal Statistical Society (B)} {\bf 46(3)}, 483-495.
\item[{\rm Barndorff-Nielsen, O.E., Cox, D. R., 1994.}] {\em Inference and Asymptotics. } Chapman and Hall, London. 
\item[{\rm Bauer, P., Hackel, P., 1978.}] The use of MOSUMS for quality control. {\em Technometrics} {\bf 20(4)}, 431-436.
\item[{\rm Bauer, P., Hackel, P., 1980.}] An extension of the MOSUM technique for quality control. {\em Technometrics} {\bf 22}, 1-7.
\item[{\rm Brazzale, A. R., Davison, A. C., Reid, N., 2006.}] {\it  Applied Asymptotics: Case Studies in Small-sample Statistics.}
 Cambridge University Press, Cambridge.
\item[{\rm Chu, C-S., Hornik, K.,  Kuan, C-M., 1995.}] MOSUM tests for parameter constancy. {\em Biometrika} {\bf 82}(3), 603-617. 
\item[{\rm Davison, A. C., Reid, N., 2022.}] The tangent exponential model.  \\
https://doi.org/10.48550/arXiv.2106.10496.
\item[{\rm Jensen, J. L., 1995.}] {\it Saddlepoint Approximations.} Oxford University Press Inc., New York.
\item[{\rm Lai, T. L., 1974.}] Control charts based on weighted sums. {\em Annals of Statistics } {\bf 2}, 134-147. 
\item[{\rm Lai, T. L., 1995.}] Sequential changepoint detection in quality control and dynamical systems (with discussions). {\em Journal of Royal Statistical Society (B)} {\bf 57(4)}, 613-658.
\item[{\rm Noonan, J., Zhigljavsky, A., 2020.}] Power of MOSUM test for online detection of transient change in mean. {\em Sequential Analysis} {\bf 39(2)}, 269-293.
\item[{\rm Severini, T. A., 2000.}] {\it Likelihood Methods in Statistics.} Oxford University Press Inc., New York.
\item [{\rm Siegmund, D., 1985.}] {\it Sequential Analysis: Tests and Confidence Intervals}. Springer, New York.
\item [{\rm Siegmund, D., 1988.}] Approximate tail probabilities for the maxima of some random fields. {\em Annals of Probability} {\bf 16(2)}, 487-501.
\item[{\rm Siegmund, D., Yakir, B., 2000.}]
Tail probabilities for the null distribution of scanning statistics. {\em Bernoulli } {\bf 6(2)}, 191-213.
\item[{\rm Siegmund, D., Yakir, B., 2008.}] Minimax optimality of the Shiryayev-Roberts change-point detection rule. {\it Journal of Statistical Planning and Inference} {\bf 138}, 2815-2825.
\item[{\rm Siegmund, D., Yakir, B., Zhang, N., 2010.}] Tail approximations for maxima of random fields by likelihood ratio transformations. {\em Sequential Analysis} {\bf 29}, 245–262.
\item[{\rm Siegmund, D., Yakir, B., Zhang, N. R., 2011.}]
Detecting simultaneous variant intervals in aligned sequences. {\em Annals of Applied Statistics} {\bf 5(2A)}, 645-668. 
\item[{\rm Tang, Y., Reid, N., 2020.}] Modified likelihood root in high dimensions. 
{\it Journal of the Royal Statistical Society Series B}  {\bf 82(5)}, 1349-1369
\item[{\rm Wu, Y., Siegmund, D., 2022.}] Sequential detection of transient signals in high dimensional data stream. Manuscript. 
\item[{\rm Xie, L., Xie, Y., Moutstakides, G. V., 2021.}] Sequential subspace change-point detection. {\em Sequential Analysis} {\bf 39(3)}, 307–335.
\item[{\rm Yakir, B., 1995.}] A note on the run length to false alarm of a change-point detection policy. {\em Annals of Statistics} {\bf 23}(1), 272-281. 
\end{description}
\end{document}